\documentclass[12pt,reqno]{amsart}

\usepackage{amssymb,enumerate,amsmath,amsthm,comment}
\usepackage{latexsym}
\usepackage{amsmath}
\usepackage{amsfonts}
\usepackage{amssymb}
\usepackage{amsthm}
\usepackage{eufrak}

\def\thebib#1{\centerline{\normalsize \sc References}
\list
{[\arabic{enumi}]}{\settowidth\labelwidth{[#1]}\leftmargin\labelwidth
 \advance\leftmargin\labelsep
  \usecounter{enumi}}
   \def\newblock{\hskip .11em plus .33em minus .07em}
    \sloppy\clubpenalty4000\widowpenalty4000
     \sfcode`\.=1000\relax}

\DeclareMathAlphabet{\E}{U}{eus}{m}{n}     

\newcommand{\PP}{{\mathbb P}}

\newcommand{\N}{{\mathbb N}}

\newcommand{\kk}{{\Bbbk}}

\newtheorem{thm}{Theorem}[section]

\theoremstyle{definition}
\newtheorem{defn}[thm]{Definition}

\newtheorem{rmk}[thm]{Remark}

\newtheorem{example}[thm]{Example}

\newtheorem{Pf}{Proof$\!\!$}         
         
\newenvironment{pf}{\begin{Pf}}{\qed\end{Pf}}

\newcommand{\rank}{\mathrm{rank}}

\DeclareMathSymbol{\twoheadrightarrow}  {\mathrel}{AMSa}{"10}

\newcounter{letter}
\renewcommand{\theletter}{\rom{(}\alph{letter}\rom{)}}

\newcounter{rnum}
\renewcommand{\thernum}{\rom{(}\roman{rnum}\rom{)}}

\begin{document}


\title[Rank of Noncommutative Quadratic Forms on Four Generators]{A Notion of Rank for 
Noncommutative\\[2mm] Quadratic Forms on Four Generators}

\baselineskip15pt

\subjclass{16S36, 15A63, 15A03}%
\keywords{quadratic form, rank, skew polynomial ring, $\mu$-symmetric
matrix%
\rule[-5mm]{0cm}{0cm}}%

\maketitle

\vspace*{0.1in}

\baselineskip15pt


\renewcommand{\thefootnote}{\fnsymbol{footnote}}
\centerline{\sc Jessica G. Cain,}
\centerline{{\sf jessicagracecain@gmail.com}}

\bigskip

\renewcommand{\thefootnote}{\fnsymbol{footnote}}
\centerline{\sc Leah R. Frauendienst,}
\centerline{Division of Mathematics and Science}
\centerline{Volunteer State Community College,
Gallatin, TN 37066}
\centerline{{\sf Leah.Frauendienst@volstate.edu}}

\bigskip
\centerline{and}
\bigskip

\centerline{\sc Padmini Veerapen}
\centerline{Department of Mathematics, ~Box 5054}
\centerline{Tennessee Tech University,
Cookeville, TN 38505}
\centerline{{\sf pveerapen@tntech.edu}}
\centerline{{\sf https://www.tntech.edu/directory/cas/math/padmini-veerapen.php}}

\bigskip
\bigskip

\begin{abstract}
\baselineskip15pt
In this paper, we extend work from \cite{VV1}, where a notion of rank, called $\mu$-rank, was proposed 
for noncommutative quadratic forms on two and three generators. In particular, we provide a definition of $\mu$-rank one and two
for noncommutative quadratic forms on four generators. 
\end{abstract}

\baselineskip21pt


\newpage

\section*{Introduction}

Quadratic forms have various applications in the commutative setting. In 2010, Cassidy and Vancliff \cite{CV1} extended the notion of a quadratic form to the noncommutative setting. They showed, in particular, that the classical correspondence between quadratic forms and symmetric matrices can be generalized to the noncommutative setting. 

Our goal in this paper is to extend work from \cite{VV1}, where a new notion of rank called, $\mu$-rank, was defined on noncommutative quadratic forms on two and three generators. In particular, we propose a definition of $\mu$-rank one and two on noncommutative quadratic forms on four generators. It was shown in \cite[Theorem 2.12]{VV2} that noncommutative quadratic forms of $\mu$-rank \textit{at most} two completely determine the point modules over graded skew Clifford algebras \cite[Definition 1.12]{CV1}. In future work, our goal is to use the definition of $\mu$-rank provided in this paper to address questions regarding graded skew Clifford algebras of global dimension four.

In Section~\ref{sec1}, we briefly review past results on $\mu$-rank of noncommutative quadratic forms on three generators. These results are used in the proofs of our main results, theorems \ref{muRk1FourVarsThm} and \ref{muRk2FourVarsThm}. In Section~\ref{sec2}, we provide proofs for theorems \ref{muRk1FourVarsThm} and \ref{muRk2FourVarsThm}, as well as a definition of $\mu$-rank one and two for noncommutative quadratic forms on four generators.

\bigskip
\bigskip


\section{Brief Summary of Work on Notion of Rank on Noncommutative Quadratic Forms}\label{sec1}


Throughout the article, $\kk$~denotes an algebraically closed field such 
that char$(\kk)\neq~2$, and $M(n,\ \kk)$ denotes the vector space of 
$n \times n$ matrices with entries in $\kk$. For a graded $\kk$-algebra~$B$, 
the span of the homogeneous elements in $B$ of degree $i$ will be denoted 
$B_i$. If $C$ is any ring or vector space, then $C^\times$ will 
denote the nonzero elements in $C$. 
\medskip

For $\{i,\ j\} \subset \{1, \ldots , n\}$, let $\mu_{ij} \in \kk^{\times}$ 
such that $\mu_{ij}\mu_{ji} = 1$ for all $i \ne j$ and we write 
$\mu = (\mu_{ij}) \in M(n,\ \kk)$.  We write $S$ for the 
quadratic $\kk$-algebra on generators $z_1, \ldots, z_n$ with defining 
relations $z_j z_i = \mu_{ij} z_i z_j$ for all $i$, $j = 1, \ldots, n$, where 
$\mu_{ii} = 1$ for all $i$. That is, 
\[S = \frac{T(V)}{\langle z_jz_i - \mu_{ij}z_iz_j : i, j = 1, \cdots, n\rangle},\]
where $T(V)$ is the tensor algebra on generators $z_1, \cdots, z_n$.

\begin{defn}\cite[\S1.2]{CV1}
\begin{enumerate}
\item[(a)]
With $\mu$ and $S$ as above, a quadratic form $Q$ is any element of $S_2$.
\item[(b)]
A matrix $M \in M(n,\ \kk)$ is called $\mu$-symmetric if $M_{ij} = 
\mu_{ij}M_{ji}$ for all $i$, $j = 1, \ldots , n$. 
\end{enumerate}
\end{defn}

We write $M^{\mu}(n, \kk)$ for the set of $\mu$-symmetric matrices in $M(n, \kk)$ and note that if $\mu_{ij} = 1$ 
for all $i, j$, then $M^{\mu}(n, \kk)$ is the set of all symmetric matrices. As shown in \cite{CV1}, the one-to-one correspondence between commutative quadratic forms and symmetric matrices can be generalized to the noncommutative setting via the one-to-one correspondence between noncommutative quadratic forms and $\mu$-symmetric matrices. In fact, $M^{\mu}(n, \kk) \cong S_2$ by mapping $M$ to $z^TMz \in S$ where $z = (z_1, ..., z_n)^T$. Theorem \ref{muRk3VarsThm} and Definition \ref{muRk3VarsDef} below reviews prior work on a definition of $\mu$-rank of quadratic forms on three generators. We also state the main results of this paper, theorems  \ref{muRk1FourVarsThm} and \ref{muRk2FourVarsThm}.

\begin{defn}\cite{VV1}\label{DefDi3Vars}
Let $M =
\left[\begin{smallmatrix}
    a & d & e \\[1mm]
    \mu_{21}d & b & f \\[1mm]
    \mu_{31}e & \mu_{32}f & c \end{smallmatrix}
\right] \in M^{\mu}(3, \kk)$ and, for $1 \le i \le 8$,
define the functions $D_i: M^{\mu}(3, \kk) \to \kk$ by\\[-3mm]
\begin{gather*}
\begin{array}{ll}
D_1(M) = 4d^2 - (1 + \mu_{12})^2 ab, \qquad &
D_4(M) = 2(1 + \mu_{23})de - (1 + \mu_{12})(1 + \mu_{13})af,\\[2mm]
D_2(M) = 4e^2 - (1 + \mu_{13})^2 ac, &
D_5(M) = 2(1 + \mu_{12})ef - (1 + \mu_{13})(1 + \mu_{23})cd,\\[2mm]
D_3(M) = 4f^2 - (1 + \mu_{23})^2 bc, &
D_6(M) = 2(1 + \mu_{13})df - (1 + \mu_{12})(1 + \mu_{23})be,
\end{array}\\[2mm]
\begin{array}{l}
D_7(M) = (\mu_{23}cd^2 - 2def + be^2)
            (\mu_{13}\mu_{21}cd^2 - 2def + \mu_{12}\mu_{23}\mu_{31}be^2),
	         \\[3mm]
D_8(M) = \mu_{21}(d + X)(e - Y) + \mu_{23}\mu_{31} (d-X)(e+Y)-2af,
\end{array}
\end{gather*}
\quad\\[-3mm]
where $X^2 = d^2 - \mu_{12} ab$ \ and \ $Y^2 = e^2 - \mu_{13} ac$.
We call $D_1, \ldots , D_6$ the $2\times 2$ $\mu$-minors of $M$ and $D_7$ and $D_8$ $\mu$-determinants of $M$.
\end{defn}

\begin{thm}\cite{VV1}\label{muRk3VarsThm}
Let $Q = az_1^2 + bz_2^2 + cz_3^2 + 2dz_1z_2 + 2ez_1z_3 + 2fz_2z_3 \in
S_2$, where $a, \ldots, f \in \kk$, and let $M \in M^{\mu}(3,\ \kk)$ be the 
$\mu$-symmetric matrix associated to $Q$. 
\begin{enumerate}[(a)]
\item[{\rm (a)}] 
There exists $L \in S_1$ such that $Q = L^2$
if and only if $D_i(M) = 0$ for all $i = 1, \ldots , 6$.
\item[{\rm (b)}] 
\begin{enumerate}
\item[{\rm (i)}] 
If $a = 0$, then there exists $L_1$, $L_2 \in S_1$ such that
$Q = L_1 L_2$ if and only if $D_7(M) = 0$;
\item[{\rm (ii)}] 
if $a \neq 0$, then there exists $L_1$, $L_2 \in S_1$ such that
$Q = L_1 L_2$ if and only if $D_8(M) = 0$ for some $X$ and $Y$ satisfying
$X^2 = d^2 - \mu_{12} ab$ \ and \ $Y^2 = e^2 - \mu_{13} ac$.
\end{enumerate}
\end{enumerate}
\end{thm}

Theorem \ref{muRk3VarsThm} leads to the following definition of $\mu$-rank for quadratic forms on three generators.

\begin{defn}\cite{VV1}\label{muRk3VarsDef}
Let $Q = az_1^2 + bz_2^2 + cz_3^2 + 2dz_1z_2 + 2ez_1z_3 + 2fz_2z_3 \in
S_2$, where $a, \ldots , f \in \kk$, with $a = 0$ or 1, 
let $M \in M^{\mu}(3,\ \kk)$ be the $\mu$-symmetric matrix associated to $Q$ 
and let $D_i : M^{\mu}(3,\ \kk) \to \kk$, for $i = 1, \ldots , 8$, be 
defined as in \ref{DefDi3Vars}.
If $n=3$, we define the function $\mu$-rank $: S_2 \to \N$ as follows:
\begin{enumerate}[(a)]
  \item if $Q = 0$, we define $\mu$-rank$(Q) = 0$;
  \item if $Q \ne 0$ and if $D_i(M) = 0$ for all $i = 1, \ldots , 6$, 
  we define $\mu$-rank$(Q) = 1$;
  \item if $D_i(M) \ne 0$ for some $i = 1, \ldots , 6$ and if 
          \[(1 - a)D_7(M) + aD_8(M) = 0,\] 
	  we define $\mu$-rank$(Q) = 2$;
  \item if $(1 - a)D_7(M) + aD_8(M) \neq 0$, we define $\mu$-rank$(Q) = 3$.
\end{enumerate}
\end{defn}

We now state the main results of this paper. Theorem \ref{muRk1FourVarsThm} allows us to define $\mu$-rank one for quadratic forms on four generators and theorem \ref{muRk2FourVarsThm} defines $\mu$-rank two. Note that the functions $D_i$, for $1 \le i \le 27$, that are referred to in theorems \ref{muRk1FourVarsThm} and \ref{muRk2FourVarsThm} are defined in \ref{DefDi4Vars}.

\newtheorem*{muRk1FourVarsThm}{Theorem \ref{muRk1FourVarsThm}}
\begin{muRk1FourVarsThm}
Let $Q = a_{11}z_1^2+a_{22}z_2^2+a_{33}z_3^2+a_{44}z_4^2+2a_{12}z_1z_2+2a_{13}z_1z_3+2a_{14}z_1z_4+2a_{23}z_2z_3+2a_{24}z_2z_4+2a_{34}z_3z_4 \in S_2$, where $a_{ij}\in\Bbbk$ for all $i$ and $j$ and let $M\in M^{\mu}(4, \Bbbk)$ be the $\mu$-symmetric matrix associated to $Q$. There exists $L\in S_{1}$ such that $Q=L^{2}$ if and only if $D_{i}(M)=0$ for all $i=1,\ldots,21$.
\end{muRk1FourVarsThm}

\newtheorem*{muRk2FourVarsThm}{Theorem \ref{muRk2FourVarsThm}}
\begin{muRk2FourVarsThm}
Let $Q = a_{11}z_1^2+a_{22}z_2^2+a_{33}z_3^2+a_{44}z_4^2+2a_{12}z_1z_2+2a_{13}z_1z_3+2a_{14}z_1z_4+2a_{23}z_2z_3+2a_{24}z_2z_4+2a_{34}z_3z_4 \in S_2$, where $a_{ij}\in\Bbbk$ for all $i$ and $j$ and let $M\in M^{\mu}(4, \Bbbk)$ be the $\mu$-symmetric matrix associated to $Q$. 
\begin{enumerate}[(a)] 
\item If $a_{11} = 0$, then there exists $L_1,L_2\in S_1$ such that $Q=L_1L_2$ if and only if $D_i(M)=0$ for all $i=22, 23, 24$. 
\item If $a_{11}\neq 0$, then there exists $L_1,L_2\in S_1$ such that $Q=L_1L_2$ if and only if $D_i(M)=0$ for all $i=25, 26, 27$ and for some $X, Y$ and $Z$ satisfying $X^2 = a_{12}^2 - \mu_{12}a_{11}a_{22}, Y^2=a_{13}^2 -\mu_{13}a_{11}a_{33}$, and $Z^2=a_{14}^2 -\mu_{14}a_{11}a_{44}$.
\end{enumerate}
\end{muRk2FourVarsThm}

Theorem \ref{muRk3VarsThm} will be useful in the proof of theorems \ref{muRk1FourVarsThm} and \ref{muRk2FourVarsThm} in section \ref{sec2}. Our goal in the next section will be to define $\mu$-rank one and two on quadratic forms on four generators. In \cite[Theorem 2.12]{VV2}, it was shown that under certain hypotheses, quadratic forms of $\mu$-rank at most two determine the number of point modules over graded skew Clifford algebras (GSCAs). The notion of rank one and two defined in section \ref{sec2} will allow the authors to apply \cite[Theorem 2.12]{VV2} to graded skew Clifford algebras of global dimension four. 

\bigskip
\bigskip
\bigskip
\bigskip

\section{$\mu$-Rank of Quadratic Forms on Four Generators}\label{sec2}

In this section, we explore further the notion of rank on noncommutative quadratic forms. Our main results are theorems \ref{muRk1FourVarsThm} and \ref{muRk2FourVarsThm}
which enable us to define $\mu$-rank one and two for quadratic forms on four generators via analogs of the determinant and minors of a $4 \times 4$ matrix.

\medskip

\begin{defn}\label{DefDi4Vars}
\noindent Let $M=
\left[\begin{smallmatrix} 
a_{11} && a_{12} && a_{13} && a_{14}\\ 
\mu_{21}a_{12} && a_{22} && a_{23} && a_{24}\\
\mu_{31}a_{13} && \mu_{32}a_{23} && a_{33} && a_{34} \\ 
\mu_{41}a_{14} && \mu_{42}a_{24} && \mu_{43}a_{34} && a_{44} 
\end{smallmatrix}\right]$
$\in M^{\mu}(4, \Bbbk)$ and, for $1\le i \le 27$, define the functions $D_{i}:M^{\mu}(4, \Bbbk) \to \Bbbk$ by\\ [-3mm]
\renewcommand*{\arraystretch}{.75}

\begin{align*}
D_1(M) &= 4(a_{12})^2 - (1+\mu_{12})^2a_{11}a_{22}, \\
D_2(M) &= 4(a_{13})^2 - (1+\mu_{13})^2a_{11}a_{33}, \\		 
D_3(M) &= 4(a_{14})^2 - (1+\mu_{14})^2a_{11}a_{44}, \\
D_4(M) &= 4(a_{23})^2 - (1+\mu_{23})^2a_{22}a_{33}, \\
D_5(M) &= 4(a_{24})^2 - (1+\mu_{24})^2a_{22}a_{44}, \\
D_6(M) &= 4(a_{34})^2 - (1+\mu_{34})^2a_{33}a_{44}, \\
D_7(M) &= 2(1+\mu_{23})a_{12}a_{13} - (1+\mu_{12})(1+\mu_{13})a_{11}a_{23}, \\
D_8(M) &= 2(1+\mu_{24})a_{12}a_{14} - (1+\mu_{12})(1+\mu_{14})a_{11}a_{24}, \\
D_9(M) &= 2(1+\mu_{13})a_{12}a_{23} - (1+\mu_{12})(1+\mu_{23})a_{13}a_{22}, \\
D_{10}(M) &= 2(1+\mu_{14})a_{12}a_{24} - (1+\mu_{12})(1+\mu_{24})a_{14}a_{22}, \end{align*} 
\begin{align*}
D_{11}(M) &= 2(1+\mu_{34})a_{13}a_{14} - (1+\mu_{13})(1+\mu_{14})a_{11}a_{34}, \\
D_{12}(M) &= 2(1+\mu_{12})a_{13}a_{23} - (1+\mu_{13})(1+\mu_{23})a_{33}a_{12}, \\
D_{13}(M) &= 2(1+\mu_{14})a_{13}a_{34} - (1+\mu_{13})(1+\mu_{34})a_{33}a_{14}, \\
D_{14}(M) &= 2(1+\mu_{12})a_{14}a_{24} - (1+\mu_{14})(1+\mu_{24})a_{12}a_{44}, \\
D_{15}(M) &= 2(1+\mu_{13})a_{14}a_{34} - (1+\mu_{14})(1+\mu_{34})a_{13}a_{44}, \\
D_{16}(M) &= 2(1+\mu_{34})a_{23}a_{24} - (1+\mu_{23})(1+\mu_{24})a_{22}a_{34}, \\
D_{17}(M) &= 2(1+\mu_{24})a_{23}a_{34} - (1+\mu_{23})(1+\mu_{34})a_{33}a_{24}, \\
D_{18}(M) &= 2(1+\mu_{23})a_{24}a_{34} - (1+\mu_{24})(1+\mu_{34})a_{23}a_{44}, \\
D_{19}(M) &= (1+\mu_{13})(1+\mu_{24})a_{12}a_{34} - (1+\mu_{12})(1+\mu_{34})a_{13}a_{24}, \\
D_{20}(M) &= (1+\mu_{14})(1+\mu_{23})a_{13}a_{24} - (1+\mu_{13})(1+\mu_{24})a_{14}a_{23}, \\
D_{21}(M) &= (1+\mu_{12})(1+\mu_{34})a_{14}a_{23} - (1+\mu_{14})(1+\mu_{23})a_{12}a_{34}, \\
D_{22}(M) &= (\mu_{23}a_{33}(a_{12})^2 -2a_{23}a_{12}a_{13} +a_{22}(a_{13})^2)(\mu_{13}\mu_{21}a_{33}(a_{12})^2 + \\
&\hspace{10ex} - 2a_{23}a_{12}a_{13} +\mu_{23}\mu_{12}\mu_{31}a_{22}(a_{13})^2), \\
D_{23}(M) &= (\mu_{24}a_{44}(a_{12})^2 -2a_{24}a_{12}a_{14} +a_{22}(a_{14})^2)(\mu_{14}\mu_{21}a_{44}(a_{12})^2 + \\
&\hspace{10ex} - 2a_{24}a_{12}a_{14} +\mu_{24}\mu_{12}\mu_{41}a_{22}(a_{14})^2), \\
D_{24}(M) &= (\mu_{34}a_{44}(a_{13})^2 -2a_{34}a_{13}a_{14} +a_{33}(a_{14})^2)(\mu_{14}\mu_{31}a_{44}(a_{13})^2 + \\
&\hspace{10ex} - 2a_{34}a_{13}a_{14} +\mu_{34}\mu_{13}\mu_{41}a_{33}(a_{14})^2), \\ 
D_{25}(M) &= \mu_{21}(a_{12} +X)(a_{13} - Y) + \mu_{23}\mu_{31}(a_{13} +Y)(a_{12} - X) - 2a_{23}a_{11}, \\
D_{26}(M) &= \mu_{21}(a_{12} +X)(a_{14} - Z) + \mu_{24}\mu_{41}(a_{14} +Z)(a_{12} - X) - 2a_{24}a_{11}, \\
D_{27}(M) &= \mu_{31}(a_{13} +Y)(a_{14} - Z) + \mu_{34}\mu_{41}(a_{14} +Z)(a_{13} - Y) - 2a_{34}a_{11}, 
\end{align*}
where $X^2 = (a_{12})^2 - \mu_{12}a_{11}a_{22}$, $Y^2=(a_{13})^2 -\mu_{13}a_{11}a_{33}$, and $Z^2=(a_{14})^2 -\mu_{14}a_{11}a_{44}$. 
\end{defn}

We call $D_{1},\ldots,D_{21}$ the $3\times 3$ $\mu$-minors of $M$ and $D_{22}, \dots, D_{27}$ the $\mu$-determinants of $M$.

\begin{thm}\label{muRk1FourVarsThm}

Let $Q = a_{11}z_1^2+a_{22}z_2^2+a_{33}z_3^2+a_{44}z_4^2+2a_{12}z_1z_2+2a_{13}z_1z_3+2a_{14}z_1z_4+2a_{23}z_2z_3+2a_{24}z_2z_4+2a_{34}z_3z_4 \in S_2$, where $a_{ij}\in\Bbbk$ for all $i$ and $j$ and let $M\in M^{\mu}(4, \Bbbk)$ be the $\mu$-symmetric matrix associated to $Q$. There exists $L\in S_{1}$ such that $Q=L^{2}$ if and only if $D_{i}(M)=0$ for all $i=1,\ldots,21$.
\begin{pf}
Suppose there exist $\alpha_{1}, \alpha_{2}, \alpha_{3}, \alpha_{4} \in \Bbbk$ such that $Q=(\alpha_{1}z_{1}+\alpha_{2}z_{2}+\alpha_{3}z_{3}+\alpha_{4}z_{4})^{2}$.  By comparing coefficients, it follows that

\begin{tabular}{rlrlrl}
(i) & $a_{11}=\alpha_{1}^{2}$, &(v)& $2a_{12}=(1+\mu_{12})\alpha_{1}\alpha_{2},$ &(viii) & $2a_{23}=(1+\mu_{23})\alpha_{2}\alpha_{3}$,\\
(ii) & $a_{22}=\alpha_{2}^{2}$, & (vi) &$2a_{13}=(1+\mu_{13})\alpha_{1}\alpha_{3}$, &(ix) &$2a_{24}=(1+\mu_{24})\alpha_{2}\alpha_{4}$,\\
(iii) & $a_{33}=\alpha_{3}^{2}$, &(vii) &$2a_{14}=(1+\mu_{14})\alpha_{1}\alpha_{4}$, &(x) &$2a_{34}=(1+\mu_{34})\alpha_{3}\alpha_{4}$, \\
(iv) &$a_{44}=\alpha_{4}^{2}$, \\
\end{tabular}\\
\vspace{-3mm}

From (v), we have
\vspace{-3mm}
\begin{equation}
4(a_{12})^{2}=(1+\mu_{12})^{2}\alpha_{1}^{2}\alpha_{2}^{2}.\nonumber 
\end{equation}
Using (i) and (ii), we have
\begin{equation}
4(a_{12})^{2}=(1+\mu_{12})^{2}a_{11}a_{22}, \nonumber
\end{equation}
so $D_{1}(M)=0$.

Similarly from (vi), we have
\vspace{-3mm}
\begin{equation}
4(a_{13})^{2}=(1+\mu_{13})^{2}\alpha_{1}^{2}\alpha_{3}^{2}.\nonumber 
\end{equation}
Using (i) and (iii), we have
\begin{equation}
4(a_{13})^{2}=(1+\mu_{13})^{2}a_{11}a_{33}, \nonumber
\end{equation}

so $D_{2}(M)=0$ and by symmetry, $D_{i}(M)=0$ for $i=3,4,5,6$. 

Using equations (v), and (vi), we have
\vspace{-3mm}
\begin{equation}
4(1+\mu_{23})a_{12}a_{13}=(1+\mu_{12})(1+\mu_{23})(1+\mu_{13})\alpha_{1}^{2}\alpha_{2}\alpha_{3}. \nonumber 
\end{equation}
Substituting in (i) and (viii), we have
\begin{equation}
4(1+\mu_{23})a_{12}a_{13}=2(1+\mu_{12})(1+\mu_{13})a_{11}a_{23}, \nonumber 
\end{equation}
so $D_{7}(M)=0$. 

Similarly using (v), (vii), we have
\vspace{-3mm}
\begin{equation}
4(1+\mu_{24})a_{12}a_{14}=(1+\mu_{12})(1+\mu_{14})(1+\mu_{24})\alpha_{1}^{2}\alpha_{2}\alpha_{4} \nonumber 
\end{equation}

Substituting (i) and (ix), we have
\begin{equation}
4(1+\mu_{24})a_{12}a_{14}=2(1+\mu_{12})(1+\mu_{14})a_{11}a_{24} \nonumber 
\end{equation}
so $D_{8}(M) = 0$ and by symmetry, $D_{i}(M)=0$ for $i=9, \ldots ,18$.

Using (v) and (x),  we have
\vspace{-3mm}
\begin{align}
4(1+\mu_{13})(1+\mu_{24})a_{12}a_{34}&=(1+\mu_{13})(1+\mu_{24})(1+\mu_{12})(1+\mu_{34})\alpha_{1}\alpha_{2}\alpha_{3}\alpha_{4}. \nonumber 
\end{align}

Substituting in (vi) and (ix), we have
\begin{align}
(1+\mu_{13})(1+\mu_{24})a_{12}a_{34}&=(1+\mu_{12})(1+\mu_{34})a_{13}a_{24}\nonumber
\end{align}

so $D_{19}(M) = 0$. 

Similarly using (vi) and (ix), we have
\vspace{-3mm}
\begin{align}
4(1+\mu_{14})(1+\mu_{23})a_{13}a_{24}&=(1+\mu_{14})(1+\mu_{23})(1+\mu_{13})(1+\mu_{24})\alpha_{1}\alpha_{2}\alpha_{3}\alpha_{4} \nonumber 
\end{align}

Substituting in (vii) and (viii), we have
\begin{align}
(1+\mu_{14})(1+\mu_{23})a_{13}a_{24}&=(1+\mu_{13})(1+\mu_{24})a_{14}a_{23} \nonumber 
\end{align}

so $D_{20}(M) = 0$ and by symmetry, $D_{21}(M)=0$. Thus, if $Q=L_{1}^{2}$ for some $L_{1}\in S_{1}$ then  $D_{i}(M)=0$ for $i=1,\ldots,21$. 


\medskip

For the converse, suppose $D_{i}(M)=0$ for all $i=1,\ldots,21$. 

\underline{Case 1.} $a_{11}=0$ 

If $a_{11}=0$, then $a_{12}=a_{13}=a_{14}=0$, since $D_{1}(M)=D_{2}(M)=D_{3}(M)=0$.  Thus, $Q$ is a quadratic form on three variables and Theorem \ref{muRk3VarsThm} applies to $Q$, that is $Q$ factors. In particular, $D_{4}(M)=D_{5}(M)=D_{6}(M)=D_{16}(M)=D_{17}(M)=D_{18}(M)=0$ and $Q = L^2$, as desired. If either $a_{22}=0$, or $a_{33}=0$, or $a_{44}=0$, the arguments follow similarly, so we omit them for brevity.




\underline{Case 2:} $a_{11}a_{22}a_{33}a_{44}\ne 0$. 

Suppose that $a_{11}a_{22}a_{33}a_{44}\ne 0$.  For simplicity, we assume that $a_{11}=1$.  Since $D_{i}(M)=0$ for all $i=1,2,3,4,5,6$, there exists $w_1, w_2, w_3, w_4, w_5, w_6\in \Bbbk$ such that

\begin{center}
\begin{tabular}{cccc}
(xi) & $2a_{12}=(1+\mu_{12})w_1$ & (xiv) & $2a_{23}=(1+\mu_{23})w_4$ \\
(xii) & $2a_{13}=(1+\mu_{13})w_2$ & (xv) & $2a_{24}=(1+\mu_{24})w_5$ \\
(xiii) & $2a_{14}=(1+\mu_{14})w_3$ & (xvi) & $2a_{34}=(1+\mu_{34})w_6$
\end{tabular}
\end{center} 
where $(w_1)^{2}=a_{22}, (w_2)^{2}=a_{33}, (w_3)^{2}=a_{44}, (w_4)^{2}=a_{22}a_{33}, (w_5)^{2}=a_{22}a_{44}, (w_6)^{2}=a_{33}a_{44}$.
Consider $Q' \in S_2$ where 
\begin{align*}
 Q'&=(z_{1}+w_1z_{2}+w_2z_{3}+w_3z_{4})^{2} \\
&=z_{1}^{2}+(w_1)^{2}z_{2}^{2}+(w_2)^{2}z_{3}^{2}+(w_3)^{2}z_{4}^{2}+w_1(1+\mu_{12})z_{1}z_{2}+w_2(1+\mu_{13})z_{1}z_{3}+w_3(1+\mu_{14})z_{1}z_{4}\\&+w_1w_2(1+\mu_{23})z_{2}z_{3}+w_1w_3(1+\mu_{24})z_{2}z_{4}+w_2w_3(1+\mu_{34})z_{3}z_{4}\\
&=z_{1}^{2}+a_{22}z_{2}^{2}+a_{33}z_{3}^{2}+a_{44}z_{4}^{2}+2a_{12}z_{1}z_{2}+2a_{13}z_{1}z_{3}+2a_{14}z_{1}z_{4}\\&+w_1w_2(1+\mu_{23})z_{2}z_{3}+w_1w_3(1+\mu_{24})z_{2}z_{4}+w_2w_3(1+\mu_{34})z_{3}z_{4}\\
\end{align*}
Since $D_{7}(M) = 0$, we have 
\begin{equation*}
2(1+\mu_{23})a_{12}a_{13}=(1+\mu_{12})(1+\mu_{13})a_{23} 
\end{equation*}
Substituting (xi), (xii), and (xiv), we have 
\begin{equation}
(1+\mu_{23})(1+\mu_{12})(1+\mu_{13})w_1w_2=(1+\mu_{12})(1+\mu_{13})(1+\mu_{23})w_4 \label{eqn1}
\end{equation}

Since $D_{8}(M) = 0$, 
\begin{equation*}
2(1+\mu_{24})a_{12}a_{14}=(1+\mu_{12})(1+\mu_{14})a_{24} 
\end{equation*}
Substituting (xi), (xiii), and (xv), we have 
\begin{equation}
(1+\mu_{24})(1+\mu_{12})(1+\mu_{14})w_1w_3=(1+\mu_{12})(1+\mu_{14})(1+\mu_{24})w_5 \label{eqn2}
\end{equation}

Since $D_{11}(M) = 0$, 
\begin{equation*}
2(1+\mu_{34})a_{13}a_{14}=(1+\mu_{13})(1+\mu_{14})a_{34} 
\end{equation*}
Substituting (xi), (xiii), and (xvi), we have 
\begin{equation}
(1+\mu_{34})(1+\mu_{12})(1+\mu_{14})w_2w_3=(1+\mu_{13})(1+\mu_{14})(1+\mu_{34})w_6 \label{eqn3}
\end{equation}

\medskip

For the rest of this proof, we split our argument into subcases where either $(1+\mu_{12})(1+\mu_{13})(1+\mu_{14})\ne0$ or $(1+\mu_{12})(1+\mu_{13})(1+\mu_{14})=0$. 

\underline{Case 2.1.} $(1+\mu_{12})(1+\mu_{13})(1+\mu_{14})\ne0$

If $(1+\mu_{12})(1+\mu_{13})(1+\mu_{14})\ne0$, then (\ref{eqn1}), (\ref{eqn2}), and (\ref{eqn3}) imply that \\

{\centering
\begin{tabular}{c}
$w_1w_2(1+\mu_{23})=2a_{23}$\\
$w_1w_3(1+\mu_{24})=2a_{24}$\\
$w_2w_3(1+\mu_{34})=2a_{34}$
\end{tabular}\\
}

Thus, $Q=Q'$.

\underline{Case 2.2.} $(1+\mu_{12})(1+\mu_{13})(1+\mu_{14})=0$  

If $(1+\mu_{12})(1+\mu_{13})(1+\mu_{14})=0$, then either $(1+\mu_{12})=0$, or $(1+\mu_{13})=0$, or $(1+\mu_{14})=0$. 

\underline{Case 2.2.1.} $(1+\mu_{12})=0$

If $(1+\mu_{12})=0$, then we may choose $w_1$ such that $w_1w_2(1+\mu_{23})=2a_{23}$.
\begin{itemize}
\item  If $(1+\mu_{13})(1+\mu_{14})\ne 0$, using $D_{11}(M) = 0$ we obtain 
\begin{align}
2(1+\mu_{34})a_{13}a_{14}&=(1+\mu_{13})(1+\mu_{14})a_{34} \nonumber \\
(1+\mu_{34})(2a_{13})(2a_{14})&=(1+\mu_{13})(1+\mu_{14})(2a_{34}) \nonumber \\
(1+\mu_{34})(1+\mu_{13})(1+\mu_{14})w_2w_3&=(1+\mu_{13})(1+\mu_{14})(2a_{34}) \nonumber \\
(1+\mu_{34})w_2w_3&=2a_{34} \nonumber 
\end{align}
Now, since $D_{20}(M) = 0$, we have
\begin{align}
(1+\mu_{14})(1+\mu_{23})a_{13}a_{24}&=(1+\mu_{13})(1+\mu_{24})a_{14}a_{23} \nonumber \\
(1+\mu_{14})(1+\mu_{23})(2a_{13})(2a_{24})&=(1+\mu_{13})(1+\mu_{24})(2a_{14})(2a_{23}) \nonumber \\
(1+\mu_{14})(1+\mu_{23})(1+\mu_{13})w_2(2a_{24})&=(1+\mu_{13})(1+\mu_{24})(1+\mu_{14})(1+\mu_{23})w_1w_2 \nonumber \\ 
(1+\mu_{23})w_2(2a_{24})&=w_1w_3(1+\mu_{24})w_2(1+\mu_{23}) \label{eqn4}
\end{align}

If $w_2(1+\mu_{23})\ne 0$, then by Equation (\ref{eqn4}), $2a_{24}=w_1w_3(1+\mu_{24})$, so that $Q'=Q$. Otherwise, if $w_2(1+\mu_{23})=0$, then either $w_2=0$ or $(1+\mu_{23})=0$. If $w_2=0$ then we may choose $w_1$ such that $2a_{24}=w_1w_3(1+\mu_{24})$, so that $Q'=Q$.  If $(1+\mu_{23})=0$ then $2a_{23}=0$. Here, instead we may choose $w_1$ such that $w_1w_3(1+\mu_{24})=2a_{24}$ and we have $Q'=Q$.

\item If $(1+\mu_{13})=0$, then we may choose $w_2$ such that $w_2w_3(1+\mu_{34})=2a_{34}$.  Since $D_{17}(M) = 0$,
\begin{align}
2(1+\mu_{24})a_{23}a_{34}&=(1+\mu_{23})(1+\mu_{34})a_{33}a_{24} \nonumber \\
(1+\mu_{24})(2a_{23})(2a_{34})&=(1+\mu_{23})(1+\mu_{34})a_{33}(2a_{24}) \nonumber \\
(1+\mu_{34})(1+\mu_{23})(1+\mu_{24})w_1w_3&=(1+\mu_{34})(1+\mu_{23})(2a_{24}) \label{eqn5}
\end{align}
If $(1+\mu_{34})(1+\mu_{23}) \ne 0$, then by Equation (\ref{eqn5}), $2a_{24}=w_1w_3(1+\mu_{24})$, so that $Q'=Q$. Otherwise, if  $(1+\mu_{34})(1+\mu_{23})=0$, then either $(1+\mu_{34})=0$ or $(1+\mu_{23})=0$ or both are zero. If $(1+\mu_{34})=0$, since $1 + \mu_{12} = 0$ then we may choose $w_1$ such that $w_1w_3(1+\mu_{24})=2a_{24}$ and $w_2$ such that $w_1w_2(1+\mu_{23})=2a_{23}$, so that $Q'=Q$.  
If $(1+\mu_{23})=0$, then $2a_{23}=0$ and we may choose $w_1$ such that $w_1w_3(1+\mu_{24})=2a_{24}$ so that $Q'=Q$, as desired.
  \\
  
\item If $(1+\mu_{14})=0$, so we may choose $w_3$ such that $w_2w_3(1+\mu_{34})=2a_{34}$ and since $D_{17}(M) = 0$, then
\begin{align}
2(1+\mu_{24})a_{23}a_{34}&=(1+\mu_{23})(1+\mu_{34})a_{33}a_{24} \nonumber \\
(1+\mu_{24})(2a_{23})(2a_{34})&=(1+\mu_{23})(1+\mu_{34})a_{33}(2a_{24}) \nonumber \\
(1+\mu_{34})(1+\mu_{23})(1+\mu_{24})w_1w_3&=(1+\mu_{34})(1+\mu_{23})(2a_{24}) \label{eqn6}
\end{align}

If $(1+\mu_{34})(1+\mu_{23}) \ne 0$, then by Equation (\ref{eqn6}), $2a_{24}=w_1w_3(1+\mu_{24})$, so that $Q'=Q$. Otherwise, either $(1+\mu_{34})=0$ or $(1+\mu_{23})=0$ or both are zero. If $(1+\mu_{34})=0$, we may chose $w_3$ such that $w_1w_3(1+\mu_{24})=2a_{24}$, so that $Q'=Q$.  
If $(1+\mu_{23})=0$, then $2a_{23}=0$ and we may choose $w_1$ such that $w_1w_3(1+\mu_{24})=2a_{24}$, so that $Q'=Q$, as desired.
\end{itemize}

{\addtolength{\leftskip}{10 mm}
\underline{Case 2..2.2.} $1+\mu_{13} = 0$ 

If $1+\mu_{13} = 0$, then we may choose $w_2$ such that $w_1w_2(1+\mu_{23})=2a_{23}$.
\begin{itemize}
 \medskip

\item 
If $(1+\mu_{12})(1+\mu_{14})\ne 0$, then since $D_{8}(M)=0$ we have,
\begin{align}
2(1+\mu_{24})a_{12}a_{14}&=(1+\mu_{12})(1+\mu_{14})a_{11}a_{24} \nonumber \\
(1+\mu_{24})(2a_{12})(2a_{14})&=(1+\mu_{12})(1+\mu_{14})(2a_{24}) \nonumber \\
(1+\mu_{12})(1+\mu_{14})(1+\mu_{24})w_1w_3&=(1+\mu_{12})(1+\mu_{14})(2a_{24}) \nonumber \\
(1+\mu_{24})w_1w_3&=2a_{24} \label{eqn7}
\end{align}
Now, since $D_{21}(M) = 0$,
\begin{align}
(1+\mu_{12})(1+\mu_{34})a_{14}a_{23}&=(1+\mu_{14})(1+\mu_{23})a_{12}a_{34} \nonumber \\
(1+\mu_{12})(1+\mu_{34})(2a_{14})(2a_{23})&=(1+\mu_{14})(1+\mu_{23})(2a_{12})(2a_{34}) \nonumber \\
(1+\mu_{12})(1+\mu_{34})(1+\mu_{14})(1+\mu_{23})w_3w_1w_2&=(1+\mu_{14})(1+\mu_{23})(1+\mu_{12})w_1(2a_{34}) \nonumber \\ 
w_1(1+\mu_{23})(1+\mu_{34})w_2w_3&=w_1(1+\mu_{23})(2a_{34}) \label{eqn8}
\end{align}

If $w_1(1+\mu_{23})\ne 0$, then by Equation (\ref{eqn8}), $2a_{34}=w_2w_3(1+\mu_{34})$ so that $Q'=Q$. Otherwise, either $w_1=0$ or $(1+\mu_{23})=0$ or both are zero. If $w_1=0$, then we may choose $w_2$ such that $2a_{34}=w_2w_3(1+\mu_{34})$ so that $Q'=Q$.  If $(1+\mu_{23})=0$, then $2a_{23}=0$. Here, instead $w_2$ could be chosen such that $w_2w_3(1+\mu_{34})=2a_{34}$ and $Q'=Q$.

\item If $(1+\mu_{12})=0$, then we may choose $w_1$ such that $w_1w_3(1+\mu_{24})=2a_{24}$ and since $D_{16}(M) = 0$, we have
\begin{align}
2(1+\mu_{34})a_{23}a_{24}&=(1+\mu_{23})(1+\mu_{24})a_{22}a_{34} \nonumber \\
(1+\mu_{34})(2a_{23})(2a_{24})&=(1+\mu_{23})(1+\mu_{24})a_{22}(2a_{34}) \nonumber \\
(1+\mu_{34})(1+\mu_{23})(1+\mu_{24})w_2w_3&=(1+\mu_{23})(1+\mu_{24})(2a_{34}) \label{eqn9}
\end{align}
If $(1+\mu_{23})(1+\mu_{24}) \ne 0$, then by Equation \ref{eqn9} $2a_{34}=w_2w_3(1+\mu_{34})$ so that $Q'=Q$. Otherwise, if $(1+\mu_{23})=0$, then $w_2$ could be chosen such that $w_2w_3(1+\mu_{34})=2a_{34}$ so that, $Q'=Q$. Similarly, if $(1+\mu_{24})=0$, then $w_2$ could be chosen such that $w_2w_3(1+\mu_{34})=2a_{34}$ and $w_1$ could have been chosen  such that $w_1w_3(1+\mu_{24})=2a_{24}$, so that $Q'=Q$.
  
\item If $(1+\mu_{14})=0$, then we may choose $w_3$ such that $w_2w_3(1+\mu_{34})=2a_{34}$ and since $D_{17}(M) = 0$ we have,
\begin{align}
2(1+\mu_{24})a_{23}a_{34}&=(1+\mu_{23})(1+\mu_{34})a_{33}a_{24} \nonumber \\
(1+\mu_{24})(2a_{23})(2a_{34})&=(1+\mu_{23})(1+\mu_{34})a_{33}(2a_{24}) \nonumber \\
(1+\mu_{34})(1+\mu_{23})(1+\mu_{24})w_1w_3&=(1+\mu_{34})(1+\mu_{23})(2a_{24}) \label{eqn10}
\end{align}

If $(1+\mu_{34})(1+\mu_{23}) \ne 0$, then by Equation (\ref{eqn10}) $2a_{24}=w_1w_3(1+\mu_{24})$ so that $Q'=Q$. Otherwise, if  $(1+\mu_{34})=0$, $w_3$ could have been chosen such that $w_1w_3(1+\mu_{24})=2a_{24}$ so that $Q'=Q$. Similarly, if $(1+\mu_{23})=0$, then $w_2$ could have been chosen such that $w_2w_3(1+\mu_{34})=2a_{34}$ and $w_3$ could have been chosen such that $w_1w_3(1+\mu_{24})=2a_{24}$ so that $Q'=Q$, as desired.
\end{itemize}
}

{\addtolength{\leftskip}{10 mm}
\underline{Case 2.2.3.} $1+\mu_{14} = 0$ 

If $1+\mu_{14} = 0$, we may choose $w_3$ such that $w_2w_3(1+\mu_{34})=2a_{34}$.
\medskip
\begin{itemize}

\item If $(1+\mu_{12})(1+\mu_{13})\ne 0$, then since $D_{7}(M)=0$ we have,
\begin{align*}
2(1+\mu_{23})a_{12}a_{13}&=(1+\mu_{12})(1+\mu_{13})a_{11}a_{23} \nonumber \\
(1+\mu_{23})(2a_{12})(2a_{13})&=(1+\mu_{12})(1+\mu_{13})(2a_{23}) \nonumber \\
(1+\mu_{12})(1+\mu_{13})(1+\mu_{23})w_1w_2&=(1+\mu_{12})(1+\mu_{13})(2a_{23}) \nonumber \\
(1+\mu_{23})w_1w_2&=2a_{23} \nonumber\\
\end{align*}
Now, since $D_{19}(M) = 0$,
\begin{align}
(1+\mu_{13})(1+\mu_{24})a_{12}a_{34}&=(1+\mu_{12})(1+\mu_{34})a_{13}a_{24} \nonumber \\
(1+\mu_{13})(1+\mu_{24})(2a_{12})(2a_{34})&=(1+\mu_{12})(1+\mu_{34})(2a_{13})(2a_{24}) \nonumber \\
(1+\mu_{12})(1+\mu_{34})(1+\mu_{13})(1+\mu_{24})w_2w_1w_3&=(1+\mu_{12})(1+\mu_{34})(1+\mu_{13})w_2(2a_{24}) \nonumber \\ 
w_2(1+\mu_{34})(1+\mu_{24})w_1w_3&=w_2(1+\mu_{34})(2a_{24}) \label{eqn11}
\end{align}

If $w_2(1+\mu_{34})\ne 0$, then by Equation (\ref{eqn11}), $2a_{24}=w_1w_3(1+\mu_{24})$ so that $Q'=Q$. Otherwise, if $w_2=0$, then $w_3$ could have been chosen such that $2a_{24}=w_1w_3(1+\mu_{24})$ so that $Q'=Q$.  If $(1+\mu_{34})=0$, $w_3$ could have been chosen  such that $w_1w_3(1+\mu_{34})=2a_{24}$ so that $Q'=Q$.

\item If $(1+\mu_{13})=0$, then we may choose $w_2$ such that $w_1w_2(1+\mu_{23})=2a_{23}$ and since $D_{17}(M) = 0$, we have
\begin{align}
2(1+\mu_{24})a_{23}a_{34}&=(1+\mu_{23})(1+\mu_{34})a_{33}a_{24} \nonumber \\
(1+\mu_{24})(2a_{23})(2a_{34})&=(1+\mu_{23})(1+\mu_{34})a_{33}(2a_{24}) \nonumber \\
(1+\mu_{23})(1+\mu_{34})(1+\mu_{24})w_1w_3&=(1+\mu_{23})(1+\mu_{34})(2a_{24}) \label{eqn12}
\end{align}
If $(1+\mu_{23})(1+\mu_{34}) \ne 0$, then by Equation (\ref{eqn12}), $2a_{24}=w_1w_3(1+\mu_{24})$ so that $Q'=Q$. Otherwise, if $(1+\mu_{23})=0$, then $w_3$ could have been chosen such that $w_1w_3(1+\mu_{24})=2a_{24}$ and $w_2$ could have been chosen such that $w_2w_3(1+\mu_{34})=2a_{34}$, so that $Q'=Q$, in both cases. Similarly, if $(1+\mu_{34})=0$, then $w_3$ could have been chosen such that $w_1w_3(1+\mu_{24})=2a_{24}$ so that, $Q'=Q$. \\
  
\item If $(1+\mu_{12})=0$, then we may choose $w_1$ such that $w_1w_2(1+\mu_{23})=2a_{23}$ and since $D_{17}(M) = 0$, we have
\begin{align}
2(1+\mu_{24})a_{23}a_{34}&=(1+\mu_{23})(1+\mu_{34})a_{33}a_{24} \nonumber \\
(1+\mu_{24})(2a_{23})(2a_{34})&=(1+\mu_{23})(1+\mu_{34})a_{33}(2a_{24}) \nonumber \\
(1+\mu_{23})(1+\mu_{34})(1+\mu_{24})w_1w_3&=(1+\mu_{23})(1+\mu_{34})(2a_{24}) \label{eqn13}
\end{align}

If $(1+\mu_{23})(1+\mu_{34}) \ne 0$, then by Equation by (\ref{eqn13}), $2a_{24}=w_1w_3(1+\mu_{24})$ so that $Q'=Q$. Otherwise, if $(1+\mu_{23})=0$, then $w_1$could have been chosen such that $w_1w_3(1+\mu_{24})=2a_{24}$ so that $Q'=Q$. Similarly, if $(1+\mu_{34})=0$, $w_3$ could have been chosen such that $w_1w_3(1+\mu_{24})=2a_{24}$ so that $Q'=Q$. 
\end{itemize}
}
Thus, when $a_{11}a_{22}a_{33}a_{44} \ne 0$ we have proved that $Q=Q'$ in each case.

Therefore, if $D_i(M) = 0$ for $1 \le i \le 21$, then $Q$ factors, as desired.
\medskip

\end{pf}
\end{thm} 


\begin{thm}\label{muRk2FourVarsThm}

Let $Q = a_{11}z_1^2+a_{22}z_2^2+a_{33}z_3^2+a_{44}z_4^2+2a_{12}z_1z_2+2a_{13}z_1z_3+2a_{14}z_1z_4+2a_{23}z_2z_3+2a_{24}z_2z_4+2a_{34}z_3z_4 \in S_2$, where $a_{ij}\in\Bbbk$ for all $i$ and $j$ and let $M\in M^{\mu}(4, \Bbbk)$ be the $\mu$-symmetric matrix associated to $Q$. 
\begin{enumerate}[(a)] 
\item If $a_{11} = 0$, then there exists $L_1,L_2\in S_1$ such that $Q=L_1L_2$ if and only if $D_i(M)=0$ for all $i=22, 23, 24$. 
\item If $a_{11}\neq 0$, then there exists $L_1,L_2\in S_1$ such that $Q=L_1L_2$ if and only if $D_i(M)=0$ for all $i=25, 26, 27$ and for some $X, Y$ and $Z$ satisfying $X^2 = a_{12}^2 - \mu_{12}a_{11}a_{22}, Y^2=a_{13}^2 -\mu_{13}a_{11}a_{33}$, and $Z^2=a_{14}^2 -\mu_{14}a_{11}a_{44}$.
\end{enumerate}

\begin{pf}
\bigskip

(a) Suppose that $a_{11} = 0$ and that the quadratic form $Q$ factors in two ways in $S$. That is, suppose there exist, $\alpha_i, \beta_j \in \Bbbk$ for $1 \le i \le 4$ and $2 \le j \le 4$ such that 
\begin{align} 
Q &= (\alpha_1z_1+ \alpha_2z_2 + \alpha_3z_3 +\alpha_4z_4)(\beta_2z_2+\beta_3z_3+\beta_4z_4) \label{eqn14} \\
&=\alpha_2\beta_2z_2^2 +\alpha_3\beta_3z_3^2+\alpha_4\beta_4z_4^2 +\alpha_1\beta_2z_1z_2+ \alpha_1\beta_3z_1z_3+\alpha_1\beta_4z_1z_4 + \nonumber\\
&+ (\alpha_2\beta_3+\mu_{23}\alpha_3\beta_2)z_2z_3 + (\alpha_2\beta_4+\mu_{24}\alpha_4\beta_2)z_2z_4 + (\alpha_3\beta_4+\mu_{34}\alpha_4\beta_3)z_3z_4.\nonumber
\end{align} By comparing coefficients, it follows that

\medskip

\begin{tabular}{ r r r r r r r }
&$\mathcal{(I)}$ &$a_{22} = \alpha_2\beta_2$, \qquad \qquad &$\mathcal{(IV)}$ &$2a_{12} = \alpha_1\beta_2$,  \qquad\qquad &$\mathcal{(VII)}$ &$2a_{23} = \alpha_2\beta_3+\mu_{23}\alpha_3\beta_2$,\\
&$\mathcal{(II)}$ &$a_{33} = \alpha_3\beta_3$, \qquad\qquad &$\mathcal{(V)}$ &$2a_{13} = \alpha_1\beta_3$, \label{case2a13} \qquad\qquad &$\mathcal{(VIII)}$ &$2a_{24} = \alpha_2\beta_4+\mu_{24}\alpha_4\beta_2$, \\
& $\mathcal{(III)}$ &$a_{44} = \alpha_4\beta_4$, \qquad\qquad &$\mathcal{(VI)}$ &$2a_{14} = \alpha_1\beta_4$, \label{case2a14} \qquad\qquad &$\mathcal{(IX)}$ &$2a_{34} = \alpha_3\beta_4+\mu_{34}\alpha_4\beta_3$. 
\end{tabular} 

\bigskip

Next, we multiply equation $\mathcal{(VII)}$ by $\mathcal{(IV)}$ and $\mathcal{(V)}$ to obtain
\begin{align*}
8a_{23}a_{12}a_{13} = (\alpha_2\beta_3+\mu_{23}\alpha_3\beta_2)(\alpha_1\beta_2)(\alpha_1\beta_3) \\
8a_{23}a_{12}a_{13} = \alpha_2\beta_3\alpha_1\beta_2\alpha_1\beta_3+\mu_{23}\alpha_3\beta_2\alpha_1\beta_2\alpha_1\beta_3.
\end{align*} We substitute in equations $\mathcal{(I)}$, $\mathcal{(II)}$ and $\mathcal{(IV)}$, $\mathcal{(V)}$ again to obtain
\begin{align*}
8a_{23}a_{12}a_{13} = a_{22}(2a_{13})^2 + \mu_{23}a_{33}(2a_{12})^2,
\end{align*} that is, 
\[\mu_{23}a_{33}a_{12}^2 -2a_{23}a_{12}a_{13} +a_{22}a_{13}^2 = 0. \]

Thus, $D_{22}(M) =0$.

Similarly, with equation $\mathcal{(VIII)}$, to obtain
\begin{align*} 
8a_{24}a_{12}a_{14} = \alpha_2\beta_4\alpha_1\beta_2\alpha_1\beta_4+\mu_{24}\alpha_4\beta_2\alpha_1\beta_2\alpha_1\beta_4 \\
8a_{24}a_{12}a_{14} = a_22(2a_{14})^2+\mu_{24}a_{44}(2a_{12})^2 \\
 \mu_{24}a_{44}a_{12}^2 -2a_{24}a_{12}a_{14} +a_{22}a_{14}^2 = 0,
\end{align*}
which implies $D_{23}(M)=0$ and with equation $\mathcal{(IX)}$ to obtain
\begin{gather*}
8a_{34}a_{13}a_{14} = \alpha_3\beta_4\alpha_1\beta_3\alpha_1\beta_4+\mu_{34}\alpha_4\beta_3\alpha_1\beta_3\alpha_1\beta_4 \\
8a_{34}a_{13}a_{14} = a_{33}(2a_{14})^2+\mu_{34}a_{44}(2a_{13})^2 \\
\mu_{34}a_{44}a_{13}^2 -2a_{34}a_{13}a_{14} +a_{33}a_{14}^2 = 0,
\end{gather*}
which implies $D_{24}(M)=0$. Hence, if $Q$ factors as in equation (\ref{eqn14}), then $D_{22}(M) = 0$, $D_{23}(M) = 0$, and $D_{24}(M) = 0$. 

\bigskip

We now consider the case when the quadratic form $Q$ factors in such a way that the order of the factors in equation (\ref{eqn14}) is switched, so suppose there exists $\alpha_i \in \Bbbk$ for $2 \le i \le 4$ and $\beta_j \in \Bbbk$ for $1 \le j \le 4$ such that
\begin{align} 
Q &= (\alpha_2z_2 + \alpha_3z_3 +\alpha_4z_4)(\beta_1z_1+\beta_2z_2+\beta_3z_3+\beta_4z_4) \label{eqn15}\\
 &= \alpha_2\beta_2z_2^2 +\alpha_3\beta_3z_3^2+\alpha_4\beta_4z_4^2 +\mu_{12}\alpha_2\beta_1z_1z_2+ \mu_{13}\alpha_3\beta_1z_1z_3+\mu_{14}\alpha_4\beta_1z_1z_4 \nonumber\\
&+ (\alpha_2\beta_3+\mu_{23}\alpha_3\beta_2)z_2z_3 + (\alpha_2\beta_4+\mu_{24}\alpha_4\beta_2)z_2z_4 + (\alpha_3\beta_4+\mu_{34}\alpha_4\beta_3)z_3z_4. \nonumber
\end{align} 
By comparing coefficients, it follows that

\medskip

\begin{tabular}{ r r r r r r r }
&$(\mathbb{I})$ & $a_{22} = \alpha_2\beta_2$, \qquad\qquad &$(\mathbb{IV})$ &$2a_{12} = \mu_{12}\alpha_2\beta_1$, \qquad\qquad &$(\mathbb{VII})$ &$2a_{23} = \alpha_2\beta_3+\mu_{23}\alpha_3\beta_2$, \label{case3eq1} \\ 
&$(\mathbb{II})$ & $a_{33} = \alpha_3\beta_3$, \qquad\qquad &$(\mathbb{V})$&$2a_{13} = \mu_{13}\alpha_3\beta_1$, \qquad\qquad &$(\mathbb{VIII})$&$2a_{24} = \alpha_2\beta_4+\mu_{24}\alpha_4\beta_2$, \label{case3eq2} \\
&$(\mathbb{III})$ &$a_{44} = \alpha_4\beta_4$, \qquad\qquad &$(\mathbb{VI})$ & $2a_{14} = \mu_{14}\alpha_4\beta_1$, \qquad\qquad &$(\mathbb{IX})$& $2a_{34} = \alpha_3\beta_4+\mu_{34}\alpha_4\beta_3$, \label{case3eq3} 
\end{tabular}

\bigskip
  
As in the first case, we multiply equation $(\mathbb{VII})$ by equations $(\mathbb{IV})$ and $(\mathbb{V})$ to obtain
\begin{align*}
8a_{23}a_{12}a_{13} &= \mu_{12}\mu_{13}\alpha_2\beta_3\alpha_2\beta_1\alpha_3\beta_1+\mu_{23}\mu_{12}\mu_{13}\alpha_3\beta_2\alpha_2\beta_1\alpha_3\beta_1
\end{align*} and substitute in equations ($\mathbb{II}$) and ($\mathbb{IV}$)
\begin{align*}
8a_{23}a_{12}a_{13} &= \mu_{12}\mu_{13}a_{33}(2\mu_{21}a_{12})^2 + \mu_{23}\mu_{12}\mu_{13} a_{22}(2\mu_{31}a_{13})^2,
\end{align*} to obtain
\begin{align*}
\mu_{13}\mu_{21}a_{33}a_{12}^2 - 2a_{23}a_{12}a_{13} +\mu_{23}\mu_{12}\mu_{31}a_{22}a_{13}^2 &= 0.
\end{align*}
Thus, $D_{22}(M) =0$. We proceed similarly to obtain
\begin{align*}
8a_{24}a_{12}a_{14} &= \mu_{12}\mu_{14}\alpha_2\beta_4\alpha_2\beta_1\alpha_4\beta_1+\mu_{24}\mu_{12}\mu_{14}\alpha_4\beta_2\alpha_2\beta_1\alpha_4\beta_1
\end{align*} that is, $\mu_{14}\mu_{21}a_{44}a_{12}^2 - 2a_{24}a_{12}a_{14} +\mu_{24}\mu_{12}\mu_{41}a_{22}a_{14}^2 = 0$, which implies that $D_{23}(M)=0$. Similarly, 
\begin{align*}
8a_{34}a_{13}a_{14} &= \mu_{13}\mu_{14}\alpha_3\beta_4\alpha_3\beta_1\alpha_4\beta_1+\mu_{34}\mu_{13}\mu_{14}\alpha_4\beta_3\alpha_3\beta_1\alpha_4\beta_1,
\end{align*} which implies that $\mu_{14}\mu_{31}a_{44}a_{13}^2 - 2a_{34}a_{13}a_{14} +\mu_{34}\mu_{13}\mu_{41}a_{33}a_{14}^2 = 0$, that is, $D_{24}(M)=0$. Therefore, if $Q$ factors as in (\ref{eqn15}), then $D_{22}(M) =0$, $D_{23}(M) =0$, and $D_{24}(M) =0$, as desired.

\medskip


For the converse, we suppose that $D_{ij}(M) = 0$ for $22 \le i, j \le 27$ with $a_{11} = 0$. 

\medskip

{\addtolength{\leftskip}{10 mm}
\underline{Case 1:} Suppose that $a_{12}a_{13}a_{14}\neq 0$. Since $D_{22}(M)=0$, either
$\mu_{23}a_{33}a_{12}^2 -2a_{23}a_{12}a_{13} +a_{22}(a_{13})^2 = 0$ or $\mu_{13}\mu_{21}a_{33}(a_{12})^2 - 2a_{23}a_{12}a_{13} +\mu_{23}\mu_{12}\mu_{31}a_{22}(a_{13})^2 = 0$. That is, either $2a_{23} =  a_{22}a_{13}a_{12}^{-1} + \mu_{23}a_{33}a_{12}a_{13}^{-1}$ or $2a_{23} = \mu_{13}\mu_{21}a_{33}a_{12}a_{13}^{-1} + \mu_{23}\mu_{12}\mu_{31}a_{22}a_{13}a_{12}^{-1}$. Now, since $D_{23}(M)=0$, then either $\mu_{24}a_{44}(a_{12})^2 -2a_{24}a_{12}a_{14} +a_{22}(a_{14})^2 = 0$ or $\mu_{14}\mu_{21}a_{44}(a_{12})^2 - 2a_{24}a_{12}a_{14} +\mu_{24}\mu_{12}\mu_{41}a_{22}(a_{14})^2 = 0$. That is, either $2a_{24} = a_{22}a_{14}a_{12}^{-1} + \mu_{24}a_{44}a_{12}a_{14}^{-1}$ or $2a_{24} =\mu_{14}\mu_{21}a_{44}a_{12}a_{14}^{-1}  +\mu_{24}\mu_{12}\mu_{41}a_{22}a_{14}a_{12}^{-1}$. Moreover, since $D_{24}(M)=0$, either $\mu_{34}a_{44}(a_{13})^2 -2a_{34}a_{13}a_{14} +a_{33}(a_{14}^2) = 0$ or $\mu_{14}\mu_{31}a_{44}(a_{13}^2) - 2a_{34}a_{13}a_{14} +\mu_{34}\mu_{13}\mu_{41}a_{33}(a_{14})^2 = 0$. That is, either $2a_{34}= a_{33}a_{14}a_{13}^{-1} + \mu_{34}a_{44}a_{13}a_{14}^{-1}$ or $2a_{34}= \mu_{14}\mu_{31}a_{44}a_{13}a_{14}^{-1} + \mu_{34}\mu_{13}\mu_{41}a_{33}a_{14}a_{13}^{-1}$. Consider the quadratic forms $Q_1,Q_2\in S_2$ such that 
\begin{align*} 
Q_1 &= (z_1+a_{22}(2a_{12})^{-1}z_2+a_{33}(2a_{13})^{-1}z_3+a_{44}(2a_{14})^{-1}z_4)(2a_{12}z_2+2a_{13}z_3+2a_{14}z_4) \\
&= a_{22}z_2^2 + a_{33}z_3^2 + a_{44}z_4^2 +2a_{12}z_1z_2+2a_{13}z_1z_3+2a_{14}z_1z_4 + (a_{22}a_{13}a_{12}^{-1} + \mu_{23}a_{33}a_{12}a_{13}^{-1})z_2z_3 \\
&\hspace{5ex}+ (a_{22}a_{14}a_{12}^{-1} + \mu_{24}a_{44}a_{12}a_{14}^{-1})z_2z_4 + (a_{33}a_{14}a_{13}^{-1} + \mu_{34}a_{44}a_{13}a_{14}^{-1})z_3z_4 \\
&= a_{22}z_2^2+a_{33}z_3^2+a_{44}z_4^2+2a_{12}z_1z_2+2a_{13}z_1z_3+2a_{14}z_1z_4+2a_{23}z_2z_3+2a_{24}z_2z_4+2a_{34}z_3z_4,
\end{align*} and 
\begin{align*}
Q_2 &= (2\mu_{21}a_{12}z_2+2\mu_{31}a_{13}z_3+2\mu_{41}a_{14}z_4)(z_1+\mu_{12}a_{22}(2a_{12})^{-1}z_2+\mu_{13}a_{33}(2a_{13})^{-1}z_3+ \\
&\hspace{5ex}+ \mu_{14}a_{44}(2a_{14})^{-1}z_4) \\
&= a_{22}z_2^2 + a_{33}z_3^2 + a_{44}z_4^2 +2a_{12}z_1z_2+2a_{13}z_1z_3+2a_{14}z_1z_4 \\
&\hspace{5ex}+ (\mu_{13}\mu_{21}a_{33}a_{12}a_{13}^{-1} + \mu_{23}\mu_{12}\mu_{31}a_{22}a_{13}a_{12}^{-1})z_2z_3 \\
&\hspace{5ex}+ (\mu_{14}\mu_{21}a_{44}a_{12}a_{14}^{-1}  +\mu_{24}\mu_{12}\mu_{41}a_{22}a_{14}a_{12}^{-1})z_2z_4 \\
&\hspace{5ex}+ (\mu_{14}\mu_{31}a_{44}a_{13}a_{14}^{-1} + \mu_{34}\mu_{13}\mu_{41}a_{33}a_{14}a_{13}^{-1})z_3z_4 \\
&= a_{22}z_2^2+a_{33}z_3^2+a_{44}z_4^2+2a_{12}z_1z_2+2a_{13}z_1z_3+2a_{14}z_1z_4+2a_{23}z_2z_3+2a_{24}z_2z_4+2a_{34}z_3z_4.
\end{align*} Thus, either $Q$ factors as in $Q_1$ or as in $Q_2$. In either case, $Q$ factors, as desired.
}

\medskip 
{\addtolength{\leftskip}{10 mm}
\underline{Case 2:} Suppose that either $a_{12}a_{13}a_{14}=0$. Without loss of generality, we consider only the following subcases since the other subcases follow by symmetry.

\begin{itemize}
\item First, we suppose that $a_{12}=0$, and $a_{13}a_{14}\neq 0$. Since $D_{23}(M) = 0$, then $\mu_{24}\mu_{12}\mu_{41}(a_{22})^2(a_{14})^4 = 0$, that is $a_{22} =0$. Now, since $D_{24}(M)=0$, then either $\mu_{34}a_{44}(a_{13})^2 -2a_{34}a_{13}a_{14} +a_{33}(a_{14})^2 = 0$ or $\mu_{14}\mu_{31}a_{44}(a_{13})^2 - 2a_{34}a_{13}a_{14} +\mu_{34}\mu_{13}\mu_{41}a_{33}(a_{14})^2 = 0$, that is, either $2a_{34}= a_{33}a_{14}a_{13}^{-1} + \mu_{34}a_{44}a_{13}a_{14}^{-1}$ or $2a_{34}= \mu_{14}\mu_{31}a_{44}a_{13}a_{14}^{-1} + \mu_{34}\mu_{13}\mu_{41}a_{33}a_{14}a_{13}^{-1}$. Since $D_{23}(M)$ and $D_{26}(M)$ are now identically zero, we may choose $a_{24} = a_{23}a_{14}a_{13}^{-1}$ or $a_{24}=\mu_{24}\mu_{32}\mu_{13}\mu_{41}a_{23}a_{14}a_{13}^{-1}$, so let $Q_1,Q_2\in S_2$ be such that 
\begin{align*} 
Q_1 &= (z_1+a_{23}a_{13}^{-1}z_2+a_{33}(2a_{13})^{-1}z_3+a_{44}(2a_{14})^{-1}z_4)(2a_{13}z_3+2a_{14}z_4) \\
&= a_{33}z_3^2 +a_{44}z_4^2 +2a_{13}z_1z_3+2a_{14}z_1z_4+2a_{23}z_2z_3+2a_{23}a_{14}a_{13}^{-1}z_2z_4 + \\
&\hspace{5ex}+ (a_{33}a_{14}a_{13}^{-1} + \mu_{34}a_{44}a_{13}a_{14}^{-1})z_3z_4 \\
&= a_{33}z_3^2 +a_{44}z_4^2 +2a_{13}z_1z_3+2a_{14}z_1z_4+2a_{23}z_2z_3+2a_{24}z_2z_4 + 2a_{34}z_3z_4
\end{align*} and 
\begin{align*} 
Q_2 &= (2\mu_{31}a_{13}z_3+2\mu_{41}a_{14}z_4)(z_1+\mu_{32}\mu_{13}a_{23}a_{13}^{-1}z_2+\mu_{13}a_{33}(2a_{13})^{-1}z_3+\mu_{14}a_{44}(2a_{14})^{-1}z_4)  \\
&= a_{33}z_3^2 +a_{44}z_4^2 +2a_{13}z_1z_3+2a_{14}z_1z_4+2a_{23}z_2z_3+2\mu_{24}\mu_{32}\mu_{13}\mu_{41}a_{23}a_{14}a_{13}^{-1}z_2z_4  \\
&\hspace{5ex}+ (\mu_{14}\mu_{31}a_{44}a_{13}a_{14}^{-1} + \mu_{34}\mu_{13}\mu_{41}a_{33}a_{14}a_{13}^{-1})z_3z_4 \\
&= a_{33}z_3^2 +a_{44}z_4^2 +2a_{13}z_1z_3+2a_{14}z_1z_4+2a_{23}z_2z_3+2a_{24}z_2z_4 + 2a_{34}z_3z_4.
\end{align*} Thus, either $Q$ factors as in $Q_1$ or as in $Q_2$. In either case, $Q$ factors, as desired. 

\item Suppose that $a_{12} = 0 = a_{13}$ and $a_{14}\neq 0$. Since $D_{23}(M)=0$, this implies that $a_{22} = 0$. Since $D_{24}(M)=0$, this implies that $a_{33} =0$. Now, consider 
\begin{align*}
Q_1 = \alpha_2\beta_2z_2^2 +\alpha_3\beta_3z_3^2+\alpha_4\beta_4z_4^2 +\alpha_1\beta_2z_1z_2+ \alpha_1\beta_3z_1z_3+\alpha_1\beta_4z_1z_4 + (\alpha_2\beta_3+\mu_{23}\alpha_3\beta_2)z_2z_3 \\
+ (\alpha_2\beta_4+\mu_{24}\alpha_4\beta_2)z_2z_4 + (\alpha_3\beta_4+\mu_{34}\alpha_4\beta_3)z_3z_4.
\end{align*}
Under our assumptions, \[\alpha_2\beta_2 = 0, \quad \alpha_3\beta_3 = 0, \quad \alpha_1\beta_2 = 0, \quad \alpha_1\beta_3 = 0, \quad \alpha_1\beta_4 \neq 0\]
which implies that $\beta_2 = 0 = \beta_3$. Hence, $\alpha_2\beta_3+\mu_{23}\alpha_3\beta_2 = 0$, that is, $a_{23}=0$. Thus, $Q = a_{44}z_4^2 +2a_{14}z_1z_4 +2a_{24}z_2z_4+2a_{34}z_3z_4 = (2a_{14}z_1+2a_{24}z_2+2a_{34}z_3 + a_{44}z_4)z_4$. 

\item Finally, suppose that $a_{12}= 0 = a_{13} = a_{14}$. The quadratic form $Q$ is now a quadratic form on three generators and the results of Section (\ref{sec1}) apply. For $Q$ to factor as a product of two linearly independent terms, we must show that $D_7(M)=0$ or $D_8(M) =0$. If $a_{22} = 0$, note that $D_{22}(M)$ corresponds to $D_7(M)$ in Definition \ref{DefDi3Vars}. Otherwise, $D_{25}(M)$ corresponds to $D_8(M)$ in Definition \ref{DefDi3Vars}. By Theorem \ref{muRk3VarsThm}, $Q$ factors, as desired. 

\end{itemize}
}


(b) Suppose that $a_{11} \neq 0$ and suppose there exist $\alpha_i, \beta_i \in \Bbbk$ for $1 \le i \le 4$  such that 
\begin{align}
Q 
&= a_{11}^{-1}(a_{11}z_1+ \alpha_2z_2 + \alpha_3z_3 +\alpha_4z_4)(a_{11}z_1+\beta_2z_2+\beta_3z_3+\beta_4z_4). \label{eqn16}
\end{align}
That is, $Q = a_{11}z_1^2 + a_{11}^{-1}\alpha_2\beta_2z_2^2 +a_{11}^{-1}\alpha_3\beta_3z_3^2+a_{11}^{-1}\alpha_4\beta_4z_4^2 +(\beta_2+\mu_{12}\alpha_2)z_1z_2+ (\beta_3+\mu_{13}\alpha_3)z_1z_3\\
+(\beta_4+\mu_{14}\alpha_4)z_1z_4 + a_{11}^{-1}(\alpha_2\beta_3+\mu_{23}\alpha_3\beta_2)z_2z_3 + a_{11}^{-1}(\alpha_2\beta_4+\mu_{24}\alpha_4\beta_2)z_2z_4 + a_{11}^{-1}(\alpha_3\beta_4+\mu_{34}\alpha_4\beta_3)z_3z_4$. 
By comparing coefficients, it follows that

\medskip

\begin{tabular}{rlrlr|}
 & (I) $a_{11}a_{22} = \alpha_2\beta_2$, \qquad \qquad &  (IV)  $2a_{12} = \beta_2+\mu_{12}\alpha_2$, \qquad \qquad & (VII) $2a_{23}a_{11} = \alpha_2\beta_3+\mu_{23}\alpha_3\beta_2$, \\
 & (II) $a_{11}a_{33} = \alpha_3\beta_3$, \qquad \qquad & (V) $2a_{13} = \beta_3+\mu_{13}\alpha_3$, \qquad \qquad & (VIII) $2a_{24}a_{11} = \alpha_2\beta_4+\mu_{24}\alpha_4\beta_2$,\\
 & (III) $a_{11}a_{44} = \alpha_4\beta_4$, \qquad \qquad & (VI) $2a_{14} = \beta_4+\mu_{14}\alpha_4$, \qquad \qquad & (IX) \qquad $2a_{34}a_{11} = \alpha_3\beta_4+\mu_{34}\alpha_4\beta_3$,
\end{tabular}\\

\medskip

From equations (IV)-(VI), we find that $\beta_2 = 2a_{12} - \mu_{12}\alpha_2$, $\beta_3 = 2a_{13} - \mu_{13}\alpha_3$, and $\beta_4 = 2a_{14} - \mu_{14}\alpha_4$. By equation (I), this implies that 
\begin{align*}
a_{11}a_{22} &= \alpha_2\beta_2 = \alpha_2(2a_{12} - \mu_{12}\alpha_2).
\end{align*} That is,
\begin{align*}
\alpha_2 &= \mu_{21}(a_{12}\pm \sqrt{a_{12}^2 -\mu_{12}a_{11}a_{22}})
\end{align*} which we write as
\begin{align*}
\alpha_2 &= \mu_{21}(a_{12}+X).
\end{align*} 
Similarly, by equations (II) \& (III), we obtain
\begin{align*}
\alpha_3 &= \mu_{31}(a_{13} +Y),
\end{align*} 
where $Y^2 = a_{13}^2-\mu_{13}a_{11}a_{33}$, and 
\begin{align*}
\alpha_4  &= \mu_{41}(a_{14} +Z).
\end{align*}
where $Z^2 = a_{14}^2 - \mu_{14}a_{11}a_{44}$.

\medskip

We may now rewrite equation (VII) solely in terms of $a_{ij}$'s, that is, we substitute $\alpha_2, \alpha_3, \beta_2, \beta_3$ in (VII) for the above expressions to obtain
\begin{align*}
\mu_{21}(a_{12} +X)(a_{13} - Y) + \mu_{23}\mu_{31}(a_{13} +Y)(a_{12} - X) - 2a_{23}a_{11} &= 0,
\end{align*}
which implies that $D_{25}(M) = 0$. Similarly, we rewrite equations (VIII) and (IX) to obtain 
\begin{align*}
\mu_{21}(a_{12} +X)(a_{14} - Z) + \mu_{24}\mu_{41}(a_{14} +Z)(a_{12} - X) - 2a_{24}a_{11} &= 0\\
\mu_{31}(a_{13} +Y)(a_{14} - Z) + \mu_{34}\mu_{41}(a_{14} +Z)(a_{13} - Y) - 2a_{34}a_{11} &= 0.
\end{align*}
Thus, $D_{26}(M) = 0$ and $D_{27}(M)=0$. Therefore, if $Q$ factors as in (\ref{eqn16}), then $D_{25}(M) = 0 = D_{26}(M) = D_{27}(M)$, as desired.  

\medskip


For the converse, we suppose that $a_{11}\neq 0$ and assume that $D_{25}(M) = 0 = D_{26}(M) = D_{27}(M)$. Since $D_{25}(M)=0$, then
\begin{align*}
2a_{23} &= a_{11}^{-1}(\mu_{21}(a_{12} +X)(a_{13} - Y) + \mu_{23}\mu_{31}(a_{13} +Y)(a_{12} - X)).
\end{align*}
Similarly, since $D_{26}(M) = 0 = D_{27}(M)$, we have that
\begin{align*}
2a_{24} &= a_{11}^{-1}(\mu_{21}(a_{12} +X)(a_{14} - Z) + \mu_{24}\mu_{41}(a_{14} +Z)(a_{12} - X))\\
2a_{34} &= a_{11}^{-1}(\mu_{31}(a_{13} +Y)(a_{14} - Z) + \mu_{34}\mu_{41}(a_{14} +Z)(a_{13} - Y)).
\end{align*}
Let $Q'\in S_2$ and let 
\begin{align*}
Q' &= a_{11}^{-1}[a_{11}z_1+ \mu_{21}(a_{12}+X)z_2 + \mu_{31}(a_{13}+Y)z_3 +\mu_{41}(a_{14} +Z)z_4]\cdot \\
&\hspace{5ex}\cdot [a_{11}z_1+(a_{12}-X)z_2+(a_{13}-Y)z_3+(a_{14}-Z)z_4] \\
&= a_{11}z_1^2 + \mu_{21}a_{11}^{-1}(a_{12}^2-X^2)z_2^2 +\mu_{31}a_{11}^{-1}(a_{13}^2-Y^2)z_3^2+\mu_{41}a_{11}^{-1}(a_{14}^2-Z^2)z_4^2 + \\
&\hspace{5ex} + 2a_{12}z_1z_2 + 2a_{13}z_1z_3+2a_{14}z_1z_4 + 2a_{23}z_2z_3 + 2a_{24}z_2z_4 + 2a_{34}z_3z_4 \\
&= a_{11}z_1^2 + a_{22}z_2^2 +a_{33}z_3^2+a_{44}z_4^2 +2a_{12}z_1z_2+ 2a_{13}z_1z_3+2a_{14}z_1z_4 + \\
&\hspace{5ex}+ 2a_{23}z_2z_3 + 2a_{24}z_2z_4 + 2a_{34}z_3z_4. 
\end{align*}
Therefore, if $a_{11} \ne 0$ and if $D_{25}(M) = 0 D_{26}(M) = D_{27}(M)$, then $Q=Q'$, as desired.
\end{pf}
\end{thm}

Theorems \ref{muRk1FourVarsThm} and \ref{muRk2FourVarsThm} suggest the following definition of $\mu$-rank for quadratic forms on four generators.

\begin{defn}\label{muRk4VarsDef}
Let $Q=a_{11}z_{1}^{2}+a_{22}z_{2}^{2}+a_{33}z_{3}^{2}+a_{44}z_{4}^{2}+2a_{12}z_{1}z_{2}+2a_{13}z_{1}z_{3}+2a_{14}z_{1}z_{4}+2a_{23}z_{2}z_{3}+2a_{24}z_{2}z_{4}+2a_{34}z_{3}z_{4}$ $\in S_{2}$, where $g_{ij}\in \Bbbk$ and let $M\in M^{\mu}(4, \Bbbk)$ be the $\mu$-symmetric matrix associated to $Q$.  Let $D_{i}:M^{\mu}(4, \Bbbk) \to \Bbbk$ for $i=1,\ldots, 27$ be defined by Definition \ref{DefDi4Vars}. We define the function $\mu$-$\rank: S_2 \to \mathbb{N}$ as follows 
\begin{enumerate}[(a)]
\item if $Q = 0$, define $\mu$-$\rank(Q)=0$;
\item if $Q\ne0$ and if $D_{i}(M)=0$ for all $i=1,\ldots,21$, define $\mu$-$\rank(Q)=1$;
\item if $D_i(M) \ne 0$ for some $i = 1, \dots, 21$ and if 
\[(1 - a_{11})D_j(M) + a_{11}D_k(M) = 0\]
for all $j = 22, 23, 24$ and for $k = 25, 26, 27$, define $\mu$-$\rank(Q)=2$ 	.
\end{enumerate}
\end{defn}

\bigskip

\begin{rmk}
When we started work on the definition of $\mu$-rank of quadratic forms on four generators, we used the \textit{Groebner} command in Maple to generate the $\mu$-minors and $\mu$-determinants defined in Definition \ref{DefDi4Vars}. However, we obtained 78 $\mu$-determinants instead of the 6 defined in our paper and were unable to apply any algorithm to reduce the 78 $\mu$-determinants to a smaller number. As a result, the method of proof in theorems \ref{muRk1FourVarsThm} and \ref{muRk2FourVarsThm} are extensions of the method used in \cite{VV1}. 
\end{rmk}

\begin{example}
Let $Q = z_1^2 + 2z_2^2 + 2z_3^2 + 2z_4^2 + 4z_1z_2 + 4z_1z_3 + 4z_1z_4 + 6z_2z_3 + 6z_2z_4 + 6z_3z_4$ and fix $\mu_{ij} = 2$ for all $1 \le i, j \le 4$. Since $D_1(M) \ne 0$, we know by Theorem \ref{muRk1FourVarsThm} that $Q$ does not factor as a perfect square. Now, since $a_{11} \ne 0$, we see that $D_{25}(M) = 0$, $D_{26}(M) = 0$, and $D_{27}(M) = 0$. By Definition \ref{muRk4VarsDef}, $\mu$-$\rank(Q) = 2$ and indeed $Q$ can be written in factored form as follows \[Q = (z_1 + z_2 + z_3 + z_4)(z_1 + 2z_2 + 2z_3 + 2z_4).\]
\end{example}

\bigskip

The notion of rank defined above could potentially extend some results in the commutative setting to the nonommutative setting. In particular, it could provide insight into the classification of quadratic Artin-Schelter regular algebras of global dimension four \cite{ATV1}.  One of the goals of finding the $\mu$-$\rank$ of quadratic forms on four generators is to extend results from \cite{SV} regarding graded Clifford algebras of global dimension four with a finite number of point modules. In \cite{VV2}, the point modules over graded skew Clifford algebras were shown to be determined by quadratic forms of $\mu$-$\rank$ at most two, using this notion of rank. Our goal now is to use the definition of $\mu$-$\rank$ provided in Definition \ref{muRk4VarsDef} to clarify the role played by quadratic forms of $\mu$-$\rank$ one and two in determining when a graded skew Clifford algebra has a finite number of point modules.  


\bigskip

\begin{thebib}{99}
\raggedbottom
\itemsep7pt
\baselineskip18pt
%

\bibitem[ATV1]{ATV1} 
{\sc M.~Artin,  J.~Tate and M.~Van den Bergh},  Some Algebras
Associated to Automorphisms of Elliptic Curves, {\it The Grothendieck
Festschrift} {\bf 1}, 33-85,
Eds.\ P.\ Cartier et al., Birkh\"auser (Boston, 1990).

\bibitem[CV1]{CV1} 
{\sc T.\ Cassidy and M.\ Vancliff}, Generalizations of Graded Clifford
Algebras and of Complete Intersections, {\em J.\ Lond.\ Math.\ Soc.}\ 
{\bf 81} (2010), 91-112.

\bibitem[SV]{SV}
{\sc D.\ R.\ Stephenson and M.\ Vancliff}, Constructing Clifford Quantum
$\PP^3$s with Finitely Many Points, {\em J.~Algebra} {\bf 312} No.~1
(2007), 86-110.

\bibitem[VV1]{VV1}{\sc M.~Vancliff and P.~P.~Veerapen}, Generalizing the Notion of Rank to Noncommutative Quadratic Forms, in ``Noncommutative Birational Geometry, Representations and Combinatorics,'' Eds. A. Berenstein and V. Retakh, {\em Contemporary Math.} {\bf 592} (2013), 241-250.

\bibitem[VV2]{VV2}
{\sc M.~Vancliff and P.~P.~Veerapen}, Point Modules over Graded Skew
Clifford Algebras, {\em J. Algebra} {\bf 420} (2014), 54-64.
\end{thebib}

\end{document}